\documentclass[11pt]{amsart}
\usepackage[T1]{fontenc}
\usepackage[utf8]{inputenc}
\usepackage[english]{babel}
\usepackage{amsmath}
\usepackage{amssymb}
\usepackage{amscd}
\begin{document}

\title{A Criterion for the Sphericity of a CR Manifold}

 \author{M. A. Stepanova}
 \address{Lomonosov Moscow State University, Faculty of Mechanics and Mathematics}
 \email{step\_masha@mail.ru}


\begin{abstract}

We formulate and prove a simple criterion for the sphericity of a CR manifold (that is, for its
equivalence to its model surface). We discuss the differences between the general case, in which
the model surface is weighted homogeneous of arbitrary degree, and the case of
Levi-nondegenerate manifolds, whose model surface is quadratic.

Bibliography: 11 titles.



\end{abstract}

\maketitle

\textbf{1. Introduction}

\vspace{3ex}

In CR geometry, a model surface is a real-algebraic CR manifold that is a natural analogue of
the tangent hyperquadric of a Levi-nondegenerate hypersurface (see Section 3 for the
definition). We formulate a simple criterion for the sphericity of a CR manifold, that is, for
its biholomorphic equivalence to its model surface. The terminology is chosen by analogy with
the case of a Levi-nondegenerate hypersurface, where sphericity means equivalence to a
nondegenerate hyperquadric in an ambient space of the same dimension. We also obtain the
following refinement of the Poincar\'e construction in CR geometry. Let $M_{0}$ be the germ at
the origin of a nondegenerate CR manifold $M$ (holomorphic nondegeneracy together with finite
Bloom--Graham type; see \cite{102}, \cite{103}), let $Q(\lambda)$ be its model surface for a
tuple of positive integer weights $\lambda$ (see Section 3), and let $Q_{0}(\lambda)$ be the
germ of $Q(\lambda)$ at the origin. The classical Poincar\'e construction yields the inequality
${\rm dim \, aut} \, M_{0}\leq{\rm dim \, aut} \, Q_{0}(\lambda)$ for the dimensions of the Lie
algebras of infinitesimal holomorphic automorphisms of the germs $M_{0}$ and $Q_{0}(\lambda)$.
We prove that $M_{0}$ is not biholomorphically equivalent to $Q_{0}$ if and only if ${\rm dim \,
aut} \, M_{0}<{\rm dim \, aut} \, Q_{0}(\lambda)$; accordingly, equality ${\rm dim \, aut} \,
M_{0}={\rm dim \, aut} \, Q_{0}(\lambda)$ holds if and only if $M_{0}$ and $Q_{0}$ are
biholomorphically equivalent. Note that the notion of sphericity depends on the choice of
weights. It is therefore appropriate to introduce the notion of \textit{$\lambda$-sphericity}:
the equivalence of the germs $M_{0}$ and $Q_{0}(\lambda)$ for some $\lambda$ (thus, a germ is
spherical if it is $\lambda$-spherical for some tuple of weights $\lambda$).
We also discuss the differences between the general case of a weighted homogeneous nondegenerate
model surface of arbitrary weight and the case of Levi-nondegenerate manifolds, whose model
surface is quadratic.

\vspace{3ex}

\textbf{2. Preliminary definitions and facts}

\vspace{3ex}

We shall need several facts from the theory of Poincar\'e--Dulac normal forms (see \cite{101}).

\textbf{Definition 1.} An ordered tuple of complex numbers
$$\lambda=(\lambda_{1},...,\lambda_{d})\in \mathbb{C}^{d}$$
is called \textit{resonant} if, for some tuple of nonnegative integers
$\alpha=(\alpha_{1},...,\alpha_{d})\in \mathbb{Z}_{+}^{d}, \,
|\alpha|=\alpha_{1}+...+\alpha_{d}\geq 2$, the resonance relation
$$\lambda_{j}=\langle \alpha,\lambda\rangle,$$
holds, where $\langle \alpha,\lambda\rangle=\alpha_{1}\lambda_{1}+...+\alpha_{d}\lambda_{d}$ is
the scalar product of the vectors $\alpha$ and $\lambda$.
A square matrix is called \textit{resonant} if the tuple of its eigenvalues, counted with
multiplicity (that is, with each eigenvalue repeated according to its multiplicity), is
resonant.
A formal vector field $F=F_{1}\frac{\partial}{\partial t_{1}}+...+F_{d}\frac{\partial}{\partial
t_{d}}$ is called \textit{resonant at the origin} if the matrix of its linear part,
$A=\Big(\frac{\partial F_{j}}{\partial t_{l}}\Big)(0), \, 1\leq j,l\leq d,$ is resonant.
The \textit{resonant vector monomial} corresponding to a resonance $\lambda_{j}=\langle
\alpha,\lambda\rangle$ is the monomial vector field $F_{j
\alpha}=t^{\alpha}\frac{\partial}{\partial t_{j}}$ (by the definition of a resonance,
$|\alpha|=\alpha_{1}+...+\alpha_{d}\geq 2$).

The \textit{Poincar\'e domain} is the set of all tuples $\lambda=(\lambda_{1},...,\lambda_{d})$
for which the origin lies outside the convex hull of $\{\lambda_{1},...,\lambda_{d}\}\subset
\mathbb{C}$.

\vspace{3ex}

In the following theorem, the notation $F=At\frac{\partial}{\partial t}+...$ denotes a
perturbation of the linear vector field $At\frac{\partial}{\partial t}$ by terms of degree two
and higher (here $t$ is the vector $(t_{1},...,t_{d}))$.

\vspace{3ex}

\textbf{Theorem 2} (Poincar\'e--Dulac theorem). Every formal vector field
$F=At\frac{\partial}{\partial t}+...$ is formally equivalent to a vector field whose linear part
has matrix equal to the Jordan normal form of $A$ and whose nonzero higher-degree terms are
resonant monomials occurring in $F$.

\vspace{3ex}

The following corollary, which will be used below, is an immediate consequence.

\vspace{3ex}

\textbf{Corollary 3.} If a vector field $F=At\frac{\partial}{\partial t}+...$ contains no
resonant monomials, then $F$ is formally equivalent to $At\frac{\partial}{\partial t}$.

\vspace{3ex}
The following theorem also holds.
\vspace{3ex}

\textbf{Theorem 4} (a special case of the Poincar\'e normalization theorem). Suppose that the
spectrum $(\lambda_{1},...,\lambda_{d})$ of the matrix $A$ belongs to the Poincar\'e domain and
that the vector field $F=At\frac{\partial}{\partial t}+...$ contains no resonant monomials. Then
$F$ is holomorphically equivalent to $At\frac{\partial}{\partial t}$.

\vspace{3ex}

\textbf{3. A criterion for sphericity}

\vspace{3ex}

Let $t=(t_{1},...,t_{d})$ be coordinates in $\mathbb{C}^{d}$.
To define a model surface, we divide the variables $t_{j}$ into two groups:
$(t_{1},...,t_{n})=(z_{1},...,z_{n})=z$ and
$(t_{n+1},...,t_{n+k})=(w_{1}=u_{1}+iv_{1},...,w_{k}=u_{k}+iv_{k})=w=u+iv$, where $n+k=d$.

Let $M\subset \mathbb{C}^{d}$ be a real-analytic manifold of CR type $(n,k)$ (that is, of CR
dimension $n$ and codimension $k$), defined in a neighborhood of the origin by the system of
equations
$$v_{j}=\rho_{j}(z,\bar{z},u), \, 1\leq j \leq k,$$
where the functions $\rho_{j}$ are real analytic.

Let $\lambda=(\lambda_{1},...,\lambda_{d})$ be a tuple of positive integers. Denote the weight
of the variable $t_{j}$ by $[t_{j}]$, and set $[t_{j}]=[\bar{t}_{j}]=\lambda_{j}, \, 1\leq j
\leq d$. This choice of weights induces a natural grading on the space of polynomials and power
series in $t$, both convergent and formal. The weight of a monomial $c \,
t_{1}^{\alpha_{1}}\cdot ... \cdot t_{d}^{\alpha_{d}}$ is $\langle \alpha,\lambda \rangle$, and a
homogeneous polynomial of weight $m$ is a sum of monomials of weight $m$.

This grading extends to vector fields with formal power series coefficients by setting
$[\frac{\partial}{\partial t_{j}}]=[\frac{\partial}{\partial \bar{t}_{j}}]=-\lambda_{j}$.

Suppose that the weights are chosen so that the defining equations of $M$ take the form
\begin{equation}\label{eq1}
v_{j}=P_{j}(z,\bar{z},u)+o(\lambda_{n+j}), \, 1\leq j \leq k,
\end{equation}
where $P_{j}(z,\bar{z},u), \, 1\leq j \leq k,$ is a polynomial of weight
$\lambda_{n+j}=[w_{j}]$, and $o(\lambda_{n+j})$ denotes terms of weight greater than
$\lambda_{n+j}$.

Under the assumption of finite Bloom--Graham type (see \cite{102}), the weights
$\lambda_{1},...,\lambda_{n}$ of the variables $z_{1},...,z_{n}$ may be assigned arbitrarily,
whereas the remaining weights $\lambda_{n+1},...,\lambda_{d}$, that is, the weights of the
variables $w_{j}, 1\leq j\leq k,$ are uniquely determined by the weights of $z_{j}, 1\leq j\leq
n$. We do not discuss in detail the procedure for choosing the weights, which is described in
\cite{110}. We merely note that, under the finite Bloom--Graham type assumption, this is a
recursive process.

Denote by
$M_{0}$ the germ of $M$ at the origin. Recall that the Lie algebra ${\rm aut} \, M_{0}$ of
infinitesimal holomorphic automorphisms of $M_{0}$ consists of vector fields tangent to $M_{0}$
of the form
$$2 \, {\rm Re} \, \Big(f_{1}(t)\frac{\partial}{\partial
t_{1}}+...+f_{d}(t)\frac{\partial}{\partial t_{d}}\Big),$$
where $f_{j}(t)$ are germs at the origin of holomorphic functions in the ambient space
$\mathbb{C}^{d}$. Such vector fields generate local one-parameter subgroups of biholomorphic
self-maps of $M_{0}$. We assume that $M_{0}$ is \textit{nondegenerate}, that is, holomorphically
nondegenerate and of finite Bloom--Graham type (see \cite{102}, \cite{103}). We do not give
definitions of these notions, since we only need the fact that they guarantee the finite
dimensionality of ${\rm aut} \, M_{0}$. We merely note that holomorphic nondegeneracy and finite
Bloom--Graham type can be verified constructively.

For the chosen tuple of weights $\lambda$, the model surface $Q(\lambda)$ of the germ $M_{0}$ is
defined by the system of equations
$$v_{j}=P_{j}(z,\bar{z},u), \, 1\leq j \leq k.$$ Clearly, the weights $\lambda_{j}$ may be
assumed to be relatively prime, since multiplying all the weights by the same factor does not
change the model surface. We make this assumption throughout.

\vspace{3ex}
Let $Q_{0}(\lambda)$ be the germ at the origin of the model surface $Q(\lambda)$. An important
role in what follows is played by the \textit{grading vector field} $E(\lambda)\in {\rm aut} \,
Q_{0}(\lambda)$, given by
$E(\lambda)=2 \, {\rm Re} \, (\lambda_{1}t_{1}\frac{\partial}{\partial
t_{1}}+...+\lambda_{d}t_{d}\frac{\partial}{\partial t_{d}})$. The grading vector field has
weight zero. It can be used to define a weight grading that coincides with the one introduced
above. Indeed, it suffices to assign to each variable $t_{j}$ the weight
$\frac{E(\lambda)(t_{j})}{t_{j}}=\lambda_{j}$. The field $E$ corresponds to the one-parameter
subgroup $\{t_{j}\longrightarrow \tau^{\lambda_{j}}t_{j}, \, 1\leq j \leq d, \, \tau>0\}$ of the
automorphism group of the model surface. The existence of this subgroup implies that every
homogeneous component of a vector field in ${\rm aut} \, Q_{0}(\lambda)$ also belongs to ${\rm
aut} \, Q_{0}(\lambda)$.

\vspace{3ex}

\textbf{Definition 5} (see \cite{110}). A nondegenerate germ $M_{0}$ is called \textit{proper}
if, for some tuple of weights $\lambda$, its model surface $Q_{0}(\lambda)$ is nondegenerate. A
germ that is not proper is called \textit{improper}.

\vspace{3ex}
Of the two nondegeneracy conditions---holomorphic nondegeneracy and finite Bloom--Graham
type---the finite-type condition is automatically inherited by the model surface, whereas
holomorphic nondegeneracy must be imposed separately.

\vspace{3ex}
Most germs are regular -- such are, for example, Levi-nondegenerate germs. An example of an improper germ is provided by the light cone in $\mathbb{C}^{3}$. This is the
tube hypersurface defined in coordinates $(\zeta_{1},\zeta_{2},\zeta_{3})$ by

$$({\rm Im} \, \zeta_{1})^{2}+({\rm Im} \, \zeta_{2})^{2}=({\rm Im} \, \zeta_{3})^{2}, \, {\rm
Im} \, \zeta_{3}>0.$$

Away from the vertex of the cone, one can choose holomorphic coordinates $(z_{1},z_{2},w=u+iv)$
in which this hypersurface is defined by a rational equation solved for $v$:

\begin{equation}\label{eq2}
v=\frac{|z_{1}|^{2}+{\rm Re} (z_{1}^{2}\bar{z}_{2})}{1-|z_{2}|^{2}}
\end{equation}
(see \cite{111}).

\vspace{3ex}






The main result of this paper is the following theorem.

\vspace{3ex}

\textbf{Theorem 6.} Let $M_{0}$ be a nondegenerate germ. Then the following conditions are
equivalent:

a) $M_{0}$ is spherical (that is, equivalent to the germ of a model surface).

b) There exists a vector field $X\in{\rm aut} \, M_{0}$ such that, in some coordinates
$(y_{1},...,y_{d})$, one has $X=2 \, {\rm Re} \, (\alpha_{1}y_{1}\frac{\partial}{\partial
y_{1}}+...+\alpha_{d}y_{d}\frac{\partial}{\partial y_{d}})$, where the $\alpha_{j}$ are positive
integers. Thus, the field is linearizable, and a linearization with positive integer eigenvalues
exists.

c) $M_{0}$ is $\lambda$-spherical (that is, holomorphically equivalent to $Q_{0}(\lambda)$ for
some tuple of weights $\lambda$).

d) ${\rm dim \, aut} \, M_{0} = {\rm dim \, aut} \, Q_{0}(\lambda)$ for some tuple of weights
$\lambda$.

e) The Lie algebras ${\rm aut} \, M_{0}$ and ${\rm aut} \, Q_{0}(\lambda)$ are isomorphic for
some tuple of weights $\lambda$.

f) ${\rm aut} \, M_{0}$ contains a field of the form $E(\lambda)+...$, where $E(\lambda)=2 \,
{\rm Re} \, (\lambda_{1}t_{1}\frac{\partial}{\partial
t_{1}}+...+\lambda_{d}t_{d}\frac{\partial}{\partial t_{d}})\in {\rm aut} \, Q_{0}(\lambda)$ is
the grading vector field of $Q_{0}(\lambda)$ for some tuple of positive integer weights
$\lambda$, and the dots in $E(\lambda)+...$ denote terms of positive weight.






\vspace{3ex}

Before proving the theorem, we give several comments and remarks (1--7) and establish some
auxiliary results.

\vspace{3ex}

1) The formulation of the sphericity criterion in part b) of the theorem is due to V. K.
Beloshapka and I. G. Kossovskiy. This condition reduces the question of the sphericity of a germ
to the linearization of a vector field, which is studied in Poincar\'e--Dulac theory (see the
facts recalled in Section 2), and depends on the presence of resonances among the eigenvalues.
In our setting, resonances are almost always present.


\vspace{3ex}

2) Part d) of Theorem 6 refines the Poincar\'e construction in CR geometry, which yields ${\rm
dim \, aut} \, M_{0}\leq {\rm dim \, aut} \, Q_{0}(\lambda)$. This estimate is not restricted to
nondegenerate germs; model surfaces for degenerate germs are defined in the same way. The
refinement is as follows: a nondegenerate germ $M_{0}$ is not biholomorphically equivalent to
$Q_{0}(\lambda)$ if and only if ${\rm dim \, aut} \, M_{0} < {\rm dim \, aut} \,
Q_{0}(\lambda)$; accordingly, equality ${\rm dim \, aut} \, M_{0}={\rm dim \, aut} \,
Q_{0}(\lambda)$ holds if and only if $M_{0}$ and $Q_{0}(\lambda)$ are biholomorphically
equivalent.

\vspace{3ex}

3) If the $\alpha_{j}$ are allowed to take zero or negative values, then $M_{0}$ may be of
infinite Bloom--Graham type, contrary to our assumption. Here are two examples.

3.1) Let $(z_{1},z_{2})$ be coordinates in $\mathbb{C}^{2}$, and let $M={\{\rm Im} \, z_{1} =
({\rm Re} \, z_{1}) |z_{2}|^{2}\}, \, \alpha_{1}=1, \alpha_{2}=0$. The germ $M_{0}$ at the
origin coincides with the germ of its model surface and is of infinite type.

3.2) Let $(z_{1},z_{2},z_{3})$ be coordinates in $\mathbb{C}^{3}$, and let $M={\{\rm Im} \,
z_{1} = ({\rm Re} \, z_{1}) {\rm Re} \, (z_{2}\bar{z}_{3})\}, \, \alpha_{1}=1, \alpha_{2}=1,
\alpha_{3}=-1$. The germ $M_{0}$ at the origin coincides with the germ of its model surface and
is of infinite type.

The Poincar\'e construction plays an important role in the proof of Theorem 6, leading to the
following natural question.

\vspace{3ex}

\textbf{Question 7.} Does the Poincar\'e construction remain valid for zero and negative
weights? More precisely, if some of the weights $\lambda_{j}$ are zero or negative, does the
inequality ${\rm dim \, aut} \, M_{0}\leq {\rm dim \, aut} \, Q_{0}(\lambda)$ hold for the germ
$Q_{0}(\lambda)$ of the model surface constructed using such a tuple of weights $\lambda$?

\vspace{3ex}

An affirmative answer to Question 7 would also suggest the following question concerning an
analogue of Theorem 6.

\vspace{3ex}

\textbf{Question 8.} Does an analogue of Theorem 6 hold for zero and negative weights?


\vspace{3ex}

Negative weights appear not yet to have found an application in CR geometry, whereas zero
weights have already proved effective (see \cite{104}, \cite{105}). Rational weights can be
reduced to integer weights, while irrational weights do not seem appropriate in this context.



\vspace{3ex}

4) One may also consider germs with singularities in their smooth structure, given by implicit
equations (for example, quadratic cones and their perturbations). In this case, additional
conditions must be imposed on the germ and its model surface to ensure that the Poincar\'e
construction is available. These conditions are as follows (see \cite{109}).

4.1) Conditions for a hypersurface:
the germ $Q_{0}(\lambda)$ is defined by a polynomial irreducible over $\mathbb{C}$;
the hypersurfaces $Q(\lambda)$ and $M$ have dimension $2d-1$ at a smooth point (that is, away
from the vertex).

4.2) Conditions for manifolds of higher codimension:
the manifolds $Q(\lambda)$ and $M$ are irreducible;
the polynomials defining $Q(\lambda)$ are irreducible over $\mathbb{C}$;
the differentials of the defining functions of $Q(\lambda)$ and $M$ are linearly independent at
a generic point (away from the vertex);
$Q(\lambda)$ and $M$ have dimension $2d-k$ at a generic point; $Q(\lambda)$ and $M$ are generic
at a generic point.

The proof of the theorem remains valid under these conditions.



\vspace{3ex}

5) In our view, parts d) and f) of the theorem are of principal interest. Moreover, the criteria
in parts d), e), and f) are fully constructive. This follows from part a) of the next lemma,
which guarantees that every germ has at most finitely many distinct model surfaces.

\vspace{3ex}

\textbf{Lemma 9.} a) A given germ $M_{0}$ has only finitely many distinct model surfaces
$Q_{0}(\lambda)$.

b) Different tuples of weights may yield the same model surface $Q_{0}(\lambda)$.

c) The number of distinct tuples of weights $\lambda$ may be infinite.


\textbf{Proof.} a) Introduce the partial order $\lesssim$ (the \textit{componentwise order}) on
$\nu$-dimensional vectors by setting $(x_{1},...,x_{\nu})\lesssim (y_{1},...,y_{\nu})$ if
$x_{j}\leq y_{j}$ for all $j$. By associating each monomial with its multidegree, this relation
induces a componentwise order on the monomials occurring in the defining equations of $M_{0}$.

We use the following fact from commutative algebra (Dickson's lemma \cite{200}, a special case
of the Hilbert basis theorem): every subset of $(\mathbb{N}\cup\{0\})^{\nu}$ has only finitely
many minimal elements with respect to the componentwise order. For any choice of weights, only
monomials minimal with respect to this order can occur in the lowest-weight component of the
defining equations. Thus, the number of possible model surfaces does not exceed the number of
subsets of the set of all minimal elements. Consequently, there are only finitely many distinct
model surfaces.

b), c) See Example 13 below.

This proves Lemma 9.

\vspace{3ex}

6) For degenerate germs, checking condition d) loses its meaning, since ${\rm dim \, aut} \,
M_{0}$ (and hence also ${\rm aut} \, Q_{0}(\lambda)$) may be infinite, even though the algebras
${\rm aut} \, M_{0}$ and ${\rm aut} \, Q_{0}(\lambda)$ may fail to be isomorphic.

\vspace{3ex}

7) A Levi-nondegenerate germ has exactly one nondegenerate model surface, which simplifies the
verification of the sphericity criterion for such manifolds.
\vspace{3ex}

We now state and prove the main auxiliary lemma.
Let $E(\lambda)=2 \, {\rm Re} \, (\lambda_{1}t_{1}\frac{\partial}{\partial
t_{1}}+...+\lambda_{d}t_{d}\frac{\partial}{\partial t_{d}})\in {\rm aut} \, Q_{0}(\lambda)$ be
the grading vector field, and let $At=\lambda_{1}t_{1}\frac{\partial}{\partial
t_{1}}+...+\lambda_{d}t_{d}\frac{\partial}{\partial t_{d}}$ be its holomorphic component, a
vector field of type $(1,0)$. Denote by $\tilde{E}(\lambda)$ a perturbation of $E(\lambda)$ by
terms of positive weight, that is, $\tilde{E}(\lambda)=E(\lambda)+...=2 \, {\rm Re}
\,(At\frac{\partial}{\partial t}+...)$, where the dots denote terms of positive weight.
\vspace{3ex}

\textbf{Lemma 10.} 1) $\tilde{E}(\lambda)$ is locally biholomorphically equivalent to its linear
part $E(\lambda)$, and

2) this biholomorphic equivalence induces a biholomorphic map between $Q_{0}(\lambda)$ and
$M_{0}$.

\textbf{Proof.} Consider the holomorphic component $At\frac{\partial}{\partial t}+...$ of the
perturbed field. The weight of each vector monomial
$c_{\alpha}t^{\alpha}\frac{\partial}{\partial t_{j}}$, where $c_{\alpha}\in\mathbb{C}$, is
$\langle \lambda,\alpha \rangle-\lambda_{j}$, and this quantity is positive for $|\alpha|\geq
2$. Thus, the perturbation contains no resonant monomials: a resonance is an equality $\langle
\lambda,\alpha \rangle-\lambda_{j}=0$ with $|\alpha|\geq 2$, so a resonant monomial must have
weight zero. This is impossible, since the perturbation consists of monomials of positive
weight.

Since all the weights $\lambda_{1},...,\lambda_{d}$ are positive, their convex hull does not
contain the origin. Hence the tuple $\lambda=(\lambda_{1},...,\lambda_{d})$ lies in the
Poincar\'e domain. By Theorem 4, the field $At\frac{\partial}{\partial t}+...$ is therefore
holomorphically equivalent to $At\frac{\partial}{\partial t}$.
It follows that $\tilde{E}(\lambda)=2 \, {\rm Re} \, (At\frac{\partial}{\partial t}+...)$ is
holomorphically equivalent to its linear part $E(\lambda)=2 \, {\rm Re} \,
(At\frac{\partial}{\partial t})$. Consequently, $M_{0}$ is holomorphically equivalent to a model
surface: the grading vector field corresponds to dilations (multiplication of all coordinates by
positive numbers), and the presence of dilations in the local automorphism group implies that
the defining equations are weighted homogeneous. Moreover, $M_{0}$ is holomorphically equivalent
to its own model surface $Q_{0}(\lambda)$, since the linear part of the corresponding
biholomorphic map does not change the model surface.


This proves Lemma 10.
\vspace{3ex}

\textbf{Remark 11.} 1) Lemma 10 can also be proved directly, without appealing to
Poincar\'e--Dulac theory. First, the formal equivalence of $E(\lambda)$ and $\tilde{E}(\lambda)$
can be verified by a direct computation, thereby establishing the formal equivalence of the
germs $M_{0}$ and $Q_{0}(\lambda)$. Next, by Theorem 1.2 of \cite{106}, every formal mapping between
holomorphically nondegenerate CR manifolds of finite type is convergent. Thus, both the fields
$E(\lambda),\tilde{E}(\lambda)$ and the germs $M_{0},Q_{0}(\lambda)$ are biholomorphically
equivalent. The fact that the tuple $\lambda$ lies in the Poincar\'e domain when the weights are
positive was pointed out by I. G. Kossovskiy.


\vspace{3ex}

\textbf{Proof of Theorem 6.} a)$\Rightarrow$b). If $M_{0}$ is equivalent to the germ of a model
surface, then in some coordinates it is defined by homogeneous polynomial equations, as ensured
by condition b).

b)$\Rightarrow$a). The grading vector field corresponds to dilations (multiplication of all
coordinates by positive numbers). The existence of dilations in the local automorphism group
implies that the defining equations are weighted homogeneous. Thus, in the coordinates
$(y_{1},...,y_{d})$, the germ is defined by homogeneous polynomial equations, which means
precisely that it is spherical.

a)$\Rightarrow$c). One can verify (see \cite{110}) that a biholomorphic equivalence of germs
induces a quasilinear equivalence of their model surfaces, that is, a biholomorphic map given by
weighted homogeneous polynomials. This yields the desired conclusion.


c)$\Rightarrow$e) and e)$\Rightarrow$d) are immediate.


d)$\Rightarrow$f). We use the following fact. If $X_{j}=\sum_{\nu=j}^{\infty}X^{(\nu)}\in {\rm
aut} \, M_{0}$, where $X^{(\nu)}$ is the component of weight $\nu$, then $X^{(j)}\in {\rm aut}
\, Q_{0}(\lambda)$ for every $\lambda$. Indeed, the lowest-weight component of the condition
that $X_{j}$ is tangent to $M_{0}$ is precisely the condition that $X^{(j)}$ is tangent to
$Q_{0}(\lambda)$. Clearly, a basis of ${\rm aut} \, M_{0}$ can be chosen so that the
lowest-weight components $X^{(j)}$ of its basis fields $X_{j}=\sum_{\nu=j}^{\infty}X^{(\nu)}\in
{\rm aut} \, M_{0}$ are linearly independent. By assumption, ${\rm dim \, aut} \, M_{0}={\rm dim
\, aut} \, Q_{0}(\lambda)$ for some $\lambda$. Hence these lowest-weight components $X^{(j)}$
form a basis of ${\rm aut} \, Q_{0}(\lambda)$. Since $E(\lambda)\in {\rm aut} \,
Q_{0}(\lambda)$, there exists a field of the form $E(\lambda)+... \in {\rm aut} \, M_{0}$.

f)$\Rightarrow$b). Apply Lemma 10. The field $E(\lambda)+...$ is transformed into the field $X$
specified in b), with $\alpha_{j}=\lambda_{j}$, and the linearizing map takes $M_{0}$ onto
$Q_{0}(\lambda)$.




This proves Theorem 6.

\vspace{3ex}
\textbf{Corollary 12.} An improper nondegenerate germ cannot be equivalent to a model surface.

\textbf{Proof.} By part c) of the theorem, $M_{0}$ can be equivalent only to one of its own
model surfaces. However, all its model surfaces are degenerate, whereas the germ itself is not.
Hence it cannot be equivalent to any model surface.

This proves Corollary 12.

\vspace{3ex}

Part d) of the theorem implies that if equality in the dimension estimate for the automorphism
algebra is attained for two model surfaces $Q_{0}(\lambda)$ and $Q_{0}(\lambda^{*})$,
corresponding to different tuples of weights $\lambda$ and $\lambda^{*}$, then $M_{0}$ is
equivalent to both. This may happen, for example, when $Q_{0}(\lambda)=Q_{0}(\lambda^{*})$, that
is, when the two germs are defined by the same equations in the same coordinate system, namely,
the one in which $M_{0}$ itself is defined. The following example shows that this is possible.

\vspace{3ex}
\textbf{Example 13.} Let $(z_{1},z_{2},w=u+iv)$ be coordinates in $\mathbb{C}^{3}$. Consider the
nondegenerate hyperquadric $\{v= 2 \, {\rm Re} \, z_{1}\bar{z}_{2}\}$. For every $\nu>1$, one
can choose the weights $\lambda_{1}=1, \lambda_{2}=\nu-1, \lambda_{3}=\nu$. For each fixed
$\nu$, these numbers are relatively prime, and there are infinitely many such tuples.
\vspace{3ex}

This raises the following question: can $M_{0}$ be simultaneously equivalent to two
\textit{distinct} model surfaces $Q_{0}(\lambda)$ and $Q_{0}(\lambda^{*})$? Here, by distinct
model surfaces we mean germs of manifolds defined in the same coordinates by different defining
relations. The answer is affirmative, as the following example shows.

\vspace{3ex}
\textbf{Example 14.} Let $(z_{1},z_{2};w=u+iv)$ be coordinates in $\mathbb{C}^{3}$, and let the
hypersurface $M$ be defined by
$$v=2 \, {\rm Re} \, (z_{1}\bar{z}_{2}+z_{2}^{2}\bar{z}_{2}).$$ Consider the two model surfaces

$$Q(\lambda)=\{v=2 \, {\rm Re} \, (z_{1}\bar{z}_{2})\}, \ \ \
\lambda=(\lambda_{1},\lambda_{2},\lambda_{3})=(1,1,2),$$

$$Q(\lambda^{*})=M, \ \ \
\lambda^{*}=(\lambda_{1}^{*},\lambda_{2}^{*},\lambda_{3}^{*})=(2,1,3).$$

The map from $Q(\lambda)$ to $Q(\lambda^{*})$ is given by $\{z_{1}=z_{1}^{*}+(z_{2}^{*})^{2}, \,
z_{2}=z_{2}^{*}, \, w=w^{*}\}$.

\vspace{3ex}
To describe, in general, the structure of a map between two model surfaces of the same germ, we
need the following definition.

\vspace{3ex}
\textbf{Definition 15.} A biholomorphic map between two germs is called \textit{quasilinear}
with respect to the weights $\lambda$ if it is given by weighted homogeneous polynomials with
respect to those weights.
\vspace{3ex}

The map between $Q_{0}(\lambda)$ and $Q_{0}(\lambda^{*})$ in Example 14 is the composition of
two different quasilinear maps with respect to $\lambda$ and $\lambda^{*}$, respectively; the
map that is quasilinear with respect to $\lambda$ is the identity. The next proposition shows
that this is the general situation.

\vspace{3ex}
\textbf{Proposition 16.} Suppose that $M_{0}$ is equivalent to both $Q_{0}(\lambda)$ and
$Q_{0}(\lambda^{*})$. Then $Q_{0}(\lambda)$ and $Q_{0}(\lambda^{*})$ are equivalent by a
composition of two quasilinear maps with respect to $\lambda$ and $\lambda^{*}$, respectively.

\textbf{Proof.}
Denote by $Q(\lambda,\lambda^{*})$ the model surface of $Q(\lambda)$ with respect to the weights
$\lambda^{*}$, and by $Q(\lambda^{*},\lambda)$ the model surface of $Q(\lambda^{*})$ with
respect to $\lambda$. A weighted homogeneous component of the defining equations is represented
by a hyperplane in the lattice $(\mathbb{N}\cup\{0\})^{2d}$, where the multidegrees
$(\alpha_{1},\bar{\alpha}_{1}, ..., \alpha_{d},\bar{\alpha}_{d})$ of the nonzero monomials
$ct_{1}^{\alpha_{1}}\bar{t}_{1}^{\bar{\alpha}_{1}}\cdot ... \cdot
t_{d}^{\alpha_{d}}\bar{t}_{d}^{\bar{\alpha}_{d}}$ occurring in the equations are marked.
Therefore, the monomials in the defining relations of $Q(\lambda,\lambda^{*})$ are specified by
the intersection of two hyperplanes. The model surface $Q(\lambda^{*},\lambda)$ is determined by
the intersection of the same hyperplanes. Hence the defining equations of
$Q(\lambda,\lambda^{*})$ and $Q(\lambda^{*},\lambda)$ coincide, so
$Q(\lambda,\lambda^{*})=Q(\lambda^{*},\lambda)$.

Next, if $M_{0}$ is equivalent to another germ $\tilde{M}_{0}$, then one can verify (see
\cite{110}) that the model surface of $M_{0}$ for the weights $\lambda$ is quasilinearly
equivalent to the model surface of $\tilde{M}_{0}$ for the same weights. Applying this statement
to $Q_{0}(\lambda)$ and $Q_{0}(\lambda^{*})$, we find that $Q_{0}(\lambda)$ is quasilinearly
equivalent to $Q_{0}(\lambda^{*},\lambda)$ with respect to $\lambda$, and that
$Q_{0}(\lambda,\lambda^{*})$ is quasilinearly equivalent to $Q_{0}(\lambda^{*})$ with respect to
$\lambda^{*}$. The desired conclusion now follows from the equality
$Q(\lambda,\lambda^{*})=Q(\lambda^{*},\lambda)$.

This proves Proposition 16.

\vspace{3ex}

The following theorem holds for Levi-nondegenerate hypersurfaces. We state it in a form
convenient for our purposes, which is weaker than the original formulation.

\vspace{3ex}


\textbf{Theorem 17} (see \cite{107}, \cite{108}). Let $M_{0}$ be the germ of a real-analytic
Levi-nondegenerate hypersurface that is not equivalent to a hyperquadric (so $M_{0}$ is not a
model germ). Then every element of the isotropy subalgebra of ${\rm aut} \, M_{0}$ is uniquely
determined by its terms of weight zero (the component $\mathfrak{g}_{0}$ of the algebra).


\vspace{3ex}

By contrast, the isotropy subalgebra of the automorphism algebra of a hyperquadric is not
determined by terms of weight zero alone: terms of weights one and two, which are nonlinear,
must also be specified.

I. G. Kossovskiy suggested that arbitrary model surfaces might admit a similar characterization.
The following example shows that a direct generalization of this fact is false.

\vspace{3ex}

\textbf{Example 18.}
Let $(z_{1},z_{2},z_{3},z_{4};w=u+iv)$ be coordinates in $\mathbb{C}^{5}$, and let $M_{0}$ be
defined by
$$\Big\{v=2 \, {\rm Re} \,
\Big(z_{1}^{2}\bar{z}_{3}+z_{2}^{2}\bar{z}_{4}\Big)+|z_{1}z_{3}|^{2}\Big\}.$$

Assign the weights $\lambda_{1}=\lambda_{2}=\lambda_{3}=\lambda_{4}=1, \, \lambda_{5}=3$. Then
$Q_{0}(\lambda)$ is defined by
$$\Big\{v=2 \, {\rm Re} \, \Big(z_{1}^{2}\bar{z}_{3}+z_{2}^{2}\bar{z}_{4}\Big)\Big\}.$$

A direct computation readily shows that ${\rm dim \, aut} \, M_{0} < {\rm dim \, aut} \,
Q_{0}(\lambda)$, so $M_{0}$ is not equivalent to its model surface. On the other hand, it is
easy to verify that the field $2 \, {\rm Re} \, \Big(iz_{2}^{2}\frac{\partial}{\partial
z_{4}}\Big)$ belongs to the isotropy subalgebra of ${\rm aut} \, M_{0}$ and has weight one. The
same field also belongs to the isotropy subalgebra of ${\rm aut} \, Q_{0}(\lambda)$. Hence the
isotropy subalgebra of this nonmodel germ $M_{0}$ is not determined by terms of weight zero.

Similar examples are readily constructed in higher codimension, for instance by taking direct
products of $M_{0}$ with nondegenerate germs.


\end{document}